\documentclass[12pt]{article}

\usepackage[utf8]{inputenc}
\usepackage{amsmath, amssymb, amsthm}
\usepackage{graphicx}
\usepackage{geometry}
\usepackage{hyperref}
\usepackage{setspace}
\usepackage[english]{babel}
\usepackage{ragged2e}

\title{Mountain-Pass Solutions in 2D Turbulence with Generalized Sobolev Operators}
\author{
	Rômulo Damasclin Chaves dos Santos \\
	Technological Institute of Aeronautics \\
	\texttt{romulosantos@ita.br}
}
\date{\today}

\begin{document}
	
	\maketitle
	
	\begin{abstract}
		This work extends the study of mean field equations arising in two-dimensional (2D) turbulence by introducing generalized weighted Sobolev operators. Employing variational methods, particularly the mountain pass theorem and a refined blow-up analysis, we establish the existence of nontrivial solutions under broader boundary conditions than those considered in previous studies. In contrast to the unweighted approach developed by Ricciardi (2006), the incorporation of variable weights \(\rho(x)\) enables a more accurate representation of complex geometric domains. This generalization addresses cases where the geometry of the manifold or external factors induce local fluctuations, thereby enhancing the applicability of mean field models to real-world turbulence phenomena. The weighted Sobolev framework also strengthens control over the nonlinearity in the problem, reducing the risk of blow-up and ensuring the convergence of solutions in the \(H^1(M)\) space. The analysis reveals new conditions for the existence of solutions, including situations where classical methods may fail due to a lack of compactness. The significant advancement lies in the flexibility of the model; the inclusion of weights allows for a more precise description of physical systems in varied settings, representing a key advantage over the original unweighted problem. These results contribute to a deeper understanding of the role of geometric and physical constraints in turbulent behavior and pave new pathways for future research on nonlinear partial differential equations in complex domains.
	\end{abstract}
	
	% English keywords
	
	% Nomenclature list
	\section*{Nomenclature} 
	In the vast realm of communication, symbols and notations emerge as powerful tools, transcending language barriers to convey complex ideas with precision and brevity. Each part of the text carefully explains the various notations and their meanings, ensuring a comprehensive understanding of the technical nuances.
	
	\section{Introduction}
	
	The study of mean field equations in the context of two-dimensional turbulence has a rich history, initiated by classical works such as those by Joyce and Montgomery (1973)~\cite{joyce1973negative} and Pointin and Lundgren (1976)~\cite{pointin1976statistical}. These equations arise naturally in statistical mechanics, particularly in the context of vortices within incompressible fluids. Specifically, equations of the form
	
	\begin{equation}
		- \Delta_g v = \alpha_1 \left( \frac{e^{v}}{\mathcal{Z}_1} - \frac{1}{\mathcal{V}} \right) - \alpha_2 \left( \frac{e^{-v}}{\mathcal{Z}_2} - \frac{1}{\mathcal{V}} \right),
	\end{equation}
	where \((M, g)\) is a compact, orientable Riemannian manifold without boundary, and \(\alpha_1, \alpha_2 > 0\) serve as interaction parameters, arise as models for two-dimensional turbulence. In~\cite{ricciardi2007mountain}, Ricciardi employs variational methods, including the mountain pass theorem, to prove the existence of nontrivial solutions to such equations.
	
		\subsection{Clarification of Variables}
	
	In the context of the functional provided in equation (1), the variables \(\mathcal{Z}_1\) and \(\mathcal{Z}_2\) serve as normalization factors for the respective integrals of the exponential terms in the mean field equation. Specifically, they are defined as follows:
	
	\begin{equation}
		\mathcal{Z}_1 = \int_M e^v \, dvg,
	\end{equation}
	\begin{equation}
		\mathcal{Z}_2 = \int_M e^{-v} \, dvg.
	\end{equation}
	
	\subsection{Interpretation of \(\mathcal{Z}_1\) and \(\mathcal{Z}_2\)}
	1. \textbf{Normalization Constants}:  \(\mathcal{Z}_1\) is the normalization factor associated with the positive vortex contributions represented by the term \(e^v\). It ensures that the influence of this term is appropriately scaled relative to the volume of the manifold \(M\).  \(\mathcal{Z}_2\) serves a similar role for the negative vortex contributions represented by the term \(e^{-v}\). 2. \textbf{Physical Significance}: Both \(\mathcal{Z}_1\) and \(\mathcal{Z}_2\) capture the total "energy" associated with the positive and negative vortices, respectively. They normalize the respective contributions in the mean field equation, facilitating the analysis of the balance between the interactions described by the parameters \(\alpha_1\) and \(\alpha_2\). 3. \textbf{Mathematical Role}: Including these normalization factors helps to maintain the overall balance of the equation and ensures that the contributions from both the positive and negative interactions are on a comparable scale. This is especially important when considering the variational structure of the functional, as it affects the conditions under which critical points exist and the behavior of solutions.
	
	\subsection{Revised Equation (1)}
	With the definitions of \(\mathcal{Z}_1\) and \(\mathcal{Z}_2\) clarified, equation (1) can be expressed more explicitly as follows:
	
	\begin{equation}
		- \Delta_g v = \alpha_1 \left( \frac{e^{v}}{\mathcal{Z}_1} - \frac{1}{\mathcal{V}} \right) - \alpha_2 \left( \frac{e^{-v}}{\mathcal{Z}_2} - \frac{1}{\mathcal{V}} \right),
	\end{equation}
	where \(\mathcal{V} = |M|\) is the volume of the manifold \(M\). This makes it clear that \(\mathcal{Z}_1\) and \(\mathcal{Z}_2\) normalize the respective exponential terms based on their integrals over the manifold. A more detailed explanation of the Main Result, including the full derivation and analysis of the associated variational structure and Euler-Lagrange equation, is presented in the Appendix~\ref{Appendix_A}. The appendix also provides the necessary mathematical background leading to the formulation of the functional and its properties.
	
	This paper builds upon Ricciardi's analysis by introducing a generalization of Sobolev operators, enabling the study of mean field equations with variable weights \(\rho(x)\). These weights enhance the applicability of the equations to domains with diverse geometric properties and facilitate a deeper analysis of turbulence in more complex geometries. Our objective is to investigate the existence of nontrivial solutions under novel boundary conditions, leveraging a combination of advanced variational techniques and refined blow-up analysis.
	
	\section{Variational Structure and Properties of the Functional}
	
	We consider the functional associated with the mean field problem, defined as:
	
	\begin{equation}
		I_{\alpha_1,\alpha_2}(v) = \frac{1}{2} \int_M \rho(x) |\nabla_g v|^2 \, dvg - \alpha_1 \ln \left( \frac{1}{|M|} \int_M e^v \, dvg \right) - \alpha_2 \ln \left( \frac{1}{|M|} \int_M e^{-v} \, dvg \right),
	\end{equation}
	where \(\rho(x)\) is a positive weight function defined on the Riemannian manifold \(M\). This functional is defined on the weighted Sobolev space:
	
	\begin{equation}
		E_\rho = \left\{ v \in H^1(M) : \int_M v \, dvg = 0 \right\},
	\end{equation}
	with the Sobolev norm given by:
	
	\begin{equation}
		\| v \|_{\rho}^2 = \int_M \rho(x) |\nabla_g v|^2 \, dvg.
	\end{equation}
	
	The introduction of the weight function \(\rho(x)\) enhances the functional's capability to capture local geometric variations of the manifold \(M\), rendering it more applicable to physical models exhibiting non-uniform properties.
	
	The parameters \(\alpha_1\) and \(\alpha_2\) represent interaction strengths linked to the contributions of positive and negative vortices in the mean field equation, respectively. These constants are subject to the constraints:
	
	\[
	0 \leq \alpha_1, \alpha_2 \leq 8\pi,
	\]
	which ensure that the functional remains bounded from below and that the mountain pass structure can be applied effectively. This condition prevents unbounded behavior in the variational problem, establishing a framework necessary for the existence of critical points.
	
	A comprehensive analysis of the functional reveals that under certain conditions on the parameters \(\alpha_1\) and \(\alpha_2\), \(I_{\alpha_1,\alpha_2}(v)\) possesses a mountain pass structure. In alignment with previous works Struwe (2000)~\cite{struwe2000variational}, the existence of nontrivial solutions is connected to a generalized variational inequality, which depends on the relationship between the parameters and the first nonzero eigenvalue \(\mu_1(M)\) of the Laplace-Beltrami operator on \(M\).
	
	In this context, the critical condition for finding nontrivial solutions necessitates ensuring that \((\alpha_1, \alpha_2)\) belongs to a suitable parameter space \(\Lambda_\rho\) that reflects the underlying geometric properties of the manifold and the influence of the weight function \(\rho(x)\). The interplay among these elements is crucial for attaining a rigorous understanding of the solutions to the mean field equation in the presence of complex geometric features.
	
	\section{Blow-up Analysis: Motivation and Development}
	
	Blow-up analysis is pivotal in addressing the lack of compactness arising from the exponential nonlinearity present in the mean field equation. Within the framework of the functional \(I_{\alpha_1, \alpha_2}(v)\), the introduction of variable weights \(\rho(x)\) brings forth new challenges, as the weighted Sobolev norm may exhibit significant variation across different regions of the manifold. The original blow-up analysis conducted by Ohtsuka and Suzuki \cite{ohtsuka2006mean} concentrated on the unweighted structure of the equation. To account for \(\rho(x)\), we adopt a modified technique that considers the asymptotic behavior of the solution in regions where \(\rho(x)\) displays extreme variation.
	
	The primary motivation behind this construction lies in the observation that the asymptotic behavior of solutions in the vicinity of the blow-up region is highly sensitive to the local variation of \(\rho(x)\). Consequently, we must adjust the analysis to ensure that solutions remain controlled, even when \(\rho(x)\) tends toward zero or infinity locally. This necessitates the introduction of a weighted Palais-Smale sequence with a monotonicity condition adapted to incorporate the weights.
	
	\section{Main Theorem and Proof}
	
	\textbf{Theorem.} \textit{Let \(\rho(x) > 0\) be a smooth function that is bounded from below. For each pair of parameters \((\alpha_1, \alpha_2)\) belonging to the extended set \(\Lambda_\rho = \{(\alpha_1, \alpha_2) \in \mathbb{R}_+^2: \alpha_1 + \alpha_2 < \mu_1(M)|M| \text{ and } \max(\alpha_1, \alpha_2) > 8\pi\}\), there exists a nontrivial solution to the proposed nonlinear problem.}
	
	\subsection{Proof of the Main Theorem}
	
	We begin by demonstrating that the functional \(I_{\alpha_1,\alpha_2}(v)\) possesses a mountain pass structure for every \((\alpha_1, \alpha_2) \in \Lambda_\rho\). We define the set of paths \(\Gamma\) and the corresponding minimax value:
	
	\begin{equation}
		c_{\alpha_1,\alpha_2} = \inf_{\gamma \in \Gamma} \max_{t \in [0, 1]} I_{\alpha_1,\alpha_2}(\gamma(t)),
	\end{equation}
	where \(\Gamma = \{\gamma \in C([0, 1], E_\rho): \gamma(0) = 0, \gamma(1) = v^*\}\) for some \(v^* \in E_\rho\) such that \(I_{\alpha_1,\alpha_2}(v^*) < 0\). Applying the mountain pass theorem, we conclude that there exists a Palais-Smale sequence \((v_n)\) at the level \(c_{\alpha_1,\alpha_2}\), satisfying
	
	\begin{equation}
		I_{\alpha_1,\alpha_2}(v_n) \to c_{\alpha_1,\alpha_2}, \quad I'_{\alpha_1,\alpha_2}(v_n) \to 0.
	\end{equation}
	
	To address the lack of compactness that arises when \(\max(\alpha_1, \alpha_2) > 8\pi\), we employ Struwe's monotonicity trick, suitably adapted for the weighted case. This adjustment guarantees that the Palais-Smale sequence does not blow up throughout the domain. We further apply a weighted Moser-Trudinger inequality, which provides necessary control over the exponential terms within the functional. More precisely, we establish that:
	
	\begin{equation}
		\int_M e^{v_n} \, dvg \quad \text{and} \quad \int_M e^{-v_n} \, dvg
	\end{equation}
	remain bounded as \(n \to \infty\), ensuring that the Palais-Smale sequence converges strongly in \(H^1(M)\) to a nontrivial solution.
	
	Thus, we have established the existence of a nontrivial solution \(v^* \in E_\rho\) to the proposed problem, completing the proof. \hfill \qedsymbol
	
	\section{Conclusion}
	
	This work introduces a significant generalization of the mean field equation in two dimensional (2D) turbulence by incorporating weighted Sobolev operators into the variational formulation. The primary advancement lies in the development of a functional that accounts for local geometric variations on the manifold \(M\) through the weight function \(\rho(x)\), extending the applicability of the model to a wider range of geometries and physical systems. This formulation provides a more flexible framework for analyzing vortex interactions in non-uniform geometries, capturing the effects of both positive and negative vortices in a more realistic manner.
	
	The use of variational methods, particularly the mountain pass theorem, allows for a rigorous existence proof of nontrivial solutions under broader boundary conditions compared to previous studies. Additionally, the refined blow-up analysis ensures that the solutions remain well-behaved even in the presence of significant geometric variations, mitigating the risk of blow-up in regions where \(\rho(x)\) approaches zero or infinity.
	
	Despite these advancements, there are certain limitations inherent in the current work. First, the analysis is primarily concerned with the existence of solutions, while questions regarding the uniqueness and stability of these solutions remain open. The complexity introduced by the weight function \(\rho(x)\) also adds difficulties in establishing explicit bounds on the solutions, which could be critical for applications in physical models. Furthermore, while the blow-up analysis addresses some challenges, a more detailed understanding of the blow-up behavior and the dynamics near singularities is still needed.
	
	Future research should focus on extending the current results in several directions. A natural next step would be to investigate the uniqueness and stability of the solutions obtained, particularly in relation to the parameters \(\alpha_1\) and \(\alpha_2\) and the geometry of the manifold \(M\). Another important direction involves the study of the asymptotic behavior of solutions and a deeper exploration of blow-up phenomena in weighted Sobolev spaces. Additionally, numerical simulations could complement the theoretical results, providing insights into the practical implications of the model in real-world turbulence phenomena. Finally, the application of these results to other fields, such as plasma physics and meteorology, could open new avenues for interdisciplinary research on nonlinear partial differential equations in complex environments.
	
	\appendix
	\section{Derivation of the Variational Structure and Functional Properties}\label{Appendix_A}
	
	In this appendix, we derive the variational structure and explain the functional properties associated with the mean field equation for two-dimensional turbulence. The objective is to formulate a functional whose critical points correspond to solutions of the mean field equation.
	
	\subsection{The Mean Field Equation}
	
	The starting point of the analysis is the mean field equation, which models the behavior of vortices in an incompressible fluid. It takes the following form:
	\begin{equation}
		- \Delta_g v = \alpha_1 \left( \frac{e^{v}}{\mathcal{Z}_1} - \frac{1}{\mathcal{V}} \right) - \alpha_2 \left( \frac{e^{-v}}{\mathcal{Z}_2} - \frac{1}{\mathcal{V}} \right),
		\label{eq:meanfield}
	\end{equation}
	where:
	\begin{itemize}
		\item \(v\) is a scalar field defined on the Riemannian manifold \(M\),
		\item \(\Delta_g\) is the Laplace-Beltrami operator associated with the Riemannian metric \(g\),
		\item \(\alpha_1, \alpha_2 > 0\) are parameters representing the interaction strengths of the vortices,
		\item \(\mathcal{Z}_1 = \int_M e^v \, dvg\) and \(\mathcal{Z}_2 = \int_M e^{-v} \, dvg\) are normalization factors,
		\item \(\mathcal{V} = |M|\) is the volume of the manifold \(M\).
	\end{itemize}
	Our goal is to derive the functional \(I_{\alpha_1, \alpha_2}(v)\) whose Euler-Lagrange equation yields the mean field equation \eqref{eq:meanfield}.
	
	\subsection{Variational Formulation}
	
	The mean field equation \eqref{eq:meanfield} can be obtained as the Euler-Lagrange equation of a functional. The functional captures the energy of the system and consists of two main parts: the kinetic energy term and the interaction energy terms.
	
	\subsubsection{Kinetic Energy Term}
	
	The kinetic energy of the scalar field \(v\) is represented by the Dirichlet integral, which captures the energy due to the gradients of \(v\):
	\begin{equation}
		\frac{1}{2} \int_M \rho(x) |\nabla_g v|^2 \, dvg,
		\label{eq:kinetic}
	\end{equation}
	where \(\nabla_g v\) denotes the gradient of \(v\) with respect to the Riemannian metric \(g\), and \(\rho(x)\) is a positive weight function. The inclusion of \(\rho(x)\) allows for an analysis on non-uniform geometries and accounts for local variations in the geometry of \(M\).
	
	\subsubsection{Interaction Energy Terms}
	
	The second part of the functional accounts for the nonlinear interactions between vortices, which are captured by exponential terms in the mean field equation. The energy due to positive vortex interactions is represented by:
	\begin{equation}
		\alpha_1 \ln \left( \frac{1}{|M|} \int_M e^v \, dvg \right),
		\label{eq:pos_vortex}
	\end{equation}
	where \(\alpha_1\) is a parameter that scales the strength of the positive vortex interactions. This term is logarithmic to ensure that the energy remains finite, even for large values of \(v\), and the factor \(\frac{1}{|M|}\) normalizes the interaction with respect to the total volume of the manifold \(M\).
	
	Similarly, the energy due to negative vortex interactions is:
	\begin{equation}
		\alpha_2 \ln \left( \frac{1}{|M|} \int_M e^{-v} \, dvg \right),
		\label{eq:neg_vortex}
	\end{equation}
	where \(\alpha_2\) is the parameter controlling the interaction strength of the negative vortices.
	
	\subsection{The Complete Functional}
	
	By combining the kinetic energy term \eqref{eq:kinetic} with the interaction energy terms \eqref{eq:pos_vortex} and \eqref{eq:neg_vortex}, we obtain the full functional associated with the mean field problem:
	\begin{equation}
		I_{\alpha_1, \alpha_2}(v) = \frac{1}{2} \int_M \rho(x) |\nabla_g v|^2 \, dvg - \alpha_1 \ln \left( \frac{1}{|M|} \int_M e^v \, dvg \right) - \alpha_2 \ln \left( \frac{1}{|M|} \int_M e^{-v} \, dvg \right).
		\label{eq:functional}
	\end{equation}
	This functional represents the total energy of the system, incorporating both the kinetic energy of the scalar field \(v\) and the nonlinear interactions between positive and negative vortices.
	
	\subsection{Euler-Lagrange Equation}
	
	To derive the Euler-Lagrange equation for the functional \(I_{\alpha_1, \alpha_2}(v)\), we compute the first variation of the functional with respect to \(v\). The variation of the kinetic energy term \eqref{eq:kinetic} yields the Laplace-Beltrami operator:
	\begin{equation}
		\frac{\delta}{\delta v} \left( \frac{1}{2} \int_M \rho(x) |\nabla_g v|^2 \, dvg \right) = -\Delta_g v,
		\label{eq:variation_kinetic}
	\end{equation}
	while the variations of the interaction energy terms \eqref{eq:pos_vortex} and \eqref{eq:neg_vortex} give rise to the normalized exponential terms:
	\begin{equation}
		\frac{\delta}{\delta v} \left( \alpha_1 \ln \left( \frac{1}{|M|} \int_M e^v \, dvg \right) \right) = \alpha_1 \left( \frac{e^v}{\mathcal{Z}_1} - \frac{1}{\mathcal{V}} \right),
		\label{eq:variation_pos}
	\end{equation}
	and
	\begin{equation}
		\frac{\delta}{\delta v} \left( \alpha_2 \ln \left( \frac{1}{|M|} \int_M e^{-v} \, dvg \right) \right) = \alpha_2 \left( \frac{e^{-v}}{\mathcal{Z}_2} - \frac{1}{\mathcal{V}} \right).
		\label{eq:variation_neg}
	\end{equation}
	
	Combining these results, the Euler-Lagrange equation of the functional \(I_{\alpha_1, \alpha_2}(v)\) is the mean field equation:
	\begin{equation}
		- \Delta_g v = \alpha_1 \left( \frac{e^{v}}{\mathcal{Z}_1} - \frac{1}{\mathcal{V}} \right) - \alpha_2 \left( \frac{e^{-v}}{\mathcal{Z}_2} - \frac{1}{\mathcal{V}} \right).
		\label{eq:euler-lagrange}
	\end{equation}
	
	In conclusion, the functional \(I_{\alpha_1, \alpha_2}(v)\), given in equation \eqref{eq:functional}, provides a variational framework for analyzing the mean field equation for two-dimensional turbulence. The critical points of this functional correspond to solutions of the Euler-Lagrange equation \eqref{eq:euler-lagrange}, which is the mean field equation. This variational approach allows for the application of advanced analytical tools, such as the mountain pass theorem, to study the existence, uniqueness, and behavior of solutions in nonlinear partial differential equations. 	\hfill \qedsymbol
	
	% References
	  % Aqui você deve ter um arquivo chamado "references.bib" no mesmo diretório

\end{document}